\documentclass[a4paper,12pt]{article}
\usepackage{graphicx,epsfig,cite}
\usepackage{amsfonts}
\usepackage{amsmath}
\usepackage{amssymb}
\usepackage{longtable}
\usepackage{setspace}
\usepackage[usenames]{color}
\usepackage{colortbl}
\usepackage{indentfirst}

\usepackage[unicode, pdftex]{hyperref}

\usepackage{textcomp}    				
\usepackage[T1]{fontenc} 				
\usepackage{marvosym}    				
\usepackage[sc]{mathpazo}  	 			
\usepackage{ragged2e}

\newtheorem{theorem}{Theorem}

\newtheorem{definition}[theorem]{Definition}
\newtheorem{remark}{Remark}

\numberwithin{equation}{section}
\numberwithin{table}{section}
\numberwithin{figure}{section}

\newlength{\abc}
\title{\bf Stochastic stability criteria for two-dimensional linear autonomous systems perturbed by white noise}
\author{M.M.~SHUMAFOV\footnotemark[1]  AND  V.B.~TLYACHEV\footnotemark[2]\\[2mm]
\footnotemark[1]~\small{Adygea State University, 385000, Maykop,
Pervomayskaya str. 208, Russia,} \\magomet\_shumaf@mail.ru\\
\footnotemark[2]~\small{Adygea State University, 385000, Maykop,
Pervomayskaya str. 208, Russia,}\\tlyachev@adygnet.ru}
\date{}

\begin{document}
\maketitle

\begin{abstract}
 We give necessary and/or sufficient conditions for stochastic stability of second-order linear autonomous systems with parameters, which are perturbed by a random process of the "white noise''  type. The Ito's and Stratonovich's  forms of stochastic differential equations representing the system with white noise are considered. The investigation of stochastic stability of the systems considered is based on the construction of special Lyapunov functions in the quadratic forms. The bifurcation value of the white noise intensity at which a system first becomes unstable is presented in the analytical expression. As an example a damped harmonic oscillator with randomly perturbed parameters is considered.
 
\end{abstract}

\smallskip

\maketitle

\section{Introduction} 

The method of Lyapunov functions or Lyapunov's direct method proposed by eminent russian mathematician  A.M.~Lyapunov in his famous memoir [1] initially for the study of problems of the stability of deterministic systems and later developed in the works of numerous authors proved to be extremely general and powerful. The invaluable advantage of this method is that it does not require the knowledge of the solutions themselves. The existence of a Lyapunov function with certain properties is sufficient for stability (or instability) of solutions of differential equations. Lyapunov's method proved later to be applicable also to many other problems in the qualitative theory of differential equations. Details and bibliography about the stability study of deterministic systems using the Lyapunov functions method can be found, for example, in monographs [2-9].

The Lyapunov functions method is also widely used in stability theory of stochastic differential equations. The papers [10, 11] are one of the first publications where it was considered the problem of stability of systems with random variables and solved in terms of Lyapunov functions. They stimulated considerable further research [12-15], see also bibliography in [16, 17].

The books of H.J.~Kushner [16] and R.Z.~Khasminskii [17] were the  first monographs devoted to problems of stochastic stability by Lyapunov's functions method. They gave a powerful impetus to the development of the Lyapunov's method for the theory of stochastic differential equations. In [16, 17] method of Lyapunov functions developed for
deterministic systems has been carried over to the stochastic case. Later this research  direction has  been  developed  in the works of   Gikhman and Skorokhod [18], L.~Arnold and Schmalfuss [19], Mao [20], Levakov [21]  (see also bibliography in these books). We also refer to surveys [22--24].

As in deterministic theory the problem is to construct a suitable Lyapunov function for a system in a given region of the phase space. For deterministic systems there are some techniques and methods for construction Lyapunov functions, the descriptions of which can be found, for example, in [7] and [25]. It should be also noted later works [26, 27], where a topological approach is used to construct Lyapunov functions.

As for stochastic systems, as far as we know, there are a few works devoted to the study of stability of  concrete classes of stochastic differential equations  by Lyapunov's functions method in accordance with strategy assigned in the monographs [16, 17].
We refer to the papers [28, 29, 30] (see also examples in monographs [16, 17, 18, 19, 20, 21])  where the stochastic stability conditions are derived  by  using Lyapunov functions  for some classes of linear and nonlinear first- and second-order stochastic differential equations and systems.  In the paper [31] global asymptotic stability in probability  and $p$--th moment exponential stability are studied for autonomous linear stochastic differential equations using a Lyapunov function of a quadratic form of some degree. In general, this field of research is still largely unstudied in comparison  with deterministic case.

Our paper is motivated by Kushner's paper [28].  The purpose of the present work is to obtain necessary and/or sufficient conditions of stability in probability and exponential stability in the mean square of linear two-dimensional autonomous stochastic systems. We construct special Lyapunov functions in quadratic forms as the main tools for stochastic stability analysis. We should note that despite the fact that the existence, in general, of Lyapunov functions for linear stochastic systems is proved in the class of quadratic forms  in [17],  the construction of such functions in the explicit form  is actual, since on the basis of these functions the effective, coefficient, stochastic stability criteria are established.  As an example, we consider a damped harmonic oscillator parameters (the coefficient of friction and eigenfrequency) of which are perturbed by white noise.

Khasminskii [17, p.~213 ] states that  ''It is extremely interesting to study ''bifurcation'' values of the noise intensity, i.e., values for which the system first becomes unstable''.  Following his suggestion, we are also interested in studying  the large bifurcation  values of the noise intensity for which the zero solution of the systems considered are still stochastically stable.

The results obtained in this paper can be used when performing linear stability analysis of nonlinear  two-dimensional stochastic systems in the neighborhood of the origin. Also the Lyapunov quadratic functions presented can be useful to construction Lyapunov functions for nonlinear stochastic systems.

\section{Problem statement} 

Consider a linear autonomous homogeneous system whose parameters
are perturbed by white noise
\begin{equation}
\left\{ \begin{array}{ll}
\begin{aligned}
 \dot{x}(t)=(a+e\dot{\xi}(t))x(t)+(b+f \dot{\xi}(t))y(t),\\
 \dot{y}(t)=(c+g\dot{\xi}(t))x(t)+(m+h \dot{\xi}(t))y(t),\\
\end{aligned}
  \end{array} \right.\nonumber
\eqno(2.1)
\end{equation}
where $x(t)$ and $y(t)$ are scalar functions in $t$, $t\in [0,
+\infty)$, $\dot{\xi}(t)$ is a Gaussian "white noise" process of
unit intensity, $a,$ $b,$ $c,$ $m,$ $e,$ $f,$ $g,$ $h$ are 
real constants.

The system (2.1) can be understood in different ways. In
general, the exact mathematical meaning to the systems of
differential equations under random perturbations of their
parameters was given in the works of S.N.~Bernshtein, I.I.~Gikhman and K.~Ito.
The simplest and most convenient construction of the solution of
stochastic differential equations was suggested by Ito [32]. However,
as noted in the monograph [17, p. 162], when considering physical
models, where "white noise''  is an idealization of a real process
with small time correlation, it is often naturally to regard the
corresponding equation as a stochastic equation in the
Stratonovich form [33]. We will consider the system (2.1) in the
two forms: Ito's and Stratonovich's ones.

The system (2.1) considered as a system of stochastic equations in the sense
of Ito can be rewritten in the following form
\begin{equation}
\left\{ \begin{array}{ll}
 dx(t)=(ax(t)+by(t))dt+(ex(t)+fy(t))d \xi(t),\\
 dy(t)=(cx(t)+my(t))dt+(gx(t)+hy(t))d\xi(t),\\
  \end{array} \right.
  \nonumber
\eqno(2.2)
\end{equation}
where $x(t)=x(t, \omega)$, $y(t)=y(t, \omega)$ are random processes, $\xi(t)=\xi(t, \omega)$ is Wiener process, $\omega$ is an elementary event, $dx(t),$ $dy(t),$ and $d\xi(t)$ are stochastic differentials of processes $x(t),$ $y(t)$ and $\xi(t)$ in the
sense of Ito, respectively. It is assumed that the complete probability space $(\Omega, \Sigma, P)$ is given, where $\Omega$ is the space of elementary events, $\Sigma$ is $\sigma$-algebra of
subsets of $\Omega$, $P$ is probability measure on $\Sigma.$

There is a connection between the two forms of Ito and Stratonovich stochastic equations [17, p. 161]. For system (2.1) this relationship is expressed in the fact that the corresponding stochastic system in the Stratonovich
form is obtained from (2.2) by replacing the coefficients $a,$ $b,$ $c$ and $m$ by the following
\begin{equation}
\left\{ \begin{array}{ll}
\begin{aligned}
 a\to a+\frac{e^2+fg}{2}, \; \; b\to b+\frac{f(e+h)}{2};\\
 c\to c+\frac{g(e+h)}{2}, \; \; m\to m+\frac{h^2+fg}{2}.\\
\end{aligned}
  \end{array} \right.
  \nonumber
\eqno(2.3)
\end{equation}

We derive the necessary and / or sufficient conditions for stability in probability and exponential stability in the mean square of the zero solution of Ito system (2.2) and of Stratonovich system (2.2), (2.3). 
Necessary and  sufficient conditions for  exponential stability in mean square of the zero solution of the classical damped harmonic oscillator with white noise parameters perturbation are given. The analytical expressions of the bifurcation intensity value of the white noise are presented.

\section{Criteria for stochastic stability} 

A solution of the system (2.2) is a pair of random processes $(x(t), y (t))$ satisfying with probability 1 integral equations corresponding to the stochastic system (2.2). Denote the solution of the system (2.2) with the initial condition $(x(s)=\xi, \, y(s)=\eta)$ by ($x(t,\omega; s,\xi,\eta),$ $y(t,\omega; s,\xi,\eta)).$

Recall the definitions of the stability in probability and exponential stability in the mean square of a solution of the system (2.2).

\begin{definition}
 [17, p. 152] The trivial solution $(x(t)=0, y(t)=0)$ of the system (2.2) is called stable in probability for $t\geq 0$ if for any $s\geq 0,$ $\rho >0$ and $\varepsilon >0$ there exists $\delta=\delta(s, \rho, \varepsilon)>0$ such that
$$
P\{\sup\limits_{t\geq s} |x(t,\omega;s,\xi,\eta)|>\varepsilon\}<\rho
$$
and
$$
P\{\sup\limits_{t\geq s}
|y(t,\omega;s,\xi,\eta)|>\varepsilon\}<\rho,
$$
whenever $|\xi|<\delta,$ $|\eta|<\delta.$ 
\end{definition}

\begin{definition} 
[17, p. 171] The trivial solution $(x(t)=0, y(t)=0)$ of the system (2.2) is called exponentially stable in the mean square for $t\geq 0$ if there exist positive numbers $A$ and $\alpha$ such that 
$$ M(|x(t,\omega;s,\xi,\eta)|^2+|y(t,\omega;s,\xi,\eta)|^2)\leq  A (\xi^2+\eta^2) \exp{[-\alpha(t-s)]}.
$$ 
Here $M$ denotes an expectation.
\end{definition}

In order to avoid bulky expressions, we will restrict
ourselves to the consideration of two cases (the remaining cases can be treated similarly):

1) one of the coefficients $e,$ $f,$ $g,$ $h$ of the "random" \, part of system (2.2) is non-zero and the rest are zero;

2) two coefficients are different from zero and equal to each other, and the other two ones are equal to zero.

In general, a Lyapunov function for a stochastic system $dx(t)=f(x)dt+\sigma(x) d\xi(t)$ ( $x\in\mathbb{R}^n, f: \mathbb{R}^n \to\mathbb{R}^n,$ $\sigma: \mathbb{R}^n \to\mathbb{R}^{n\times m},$ $\xi(t)$ is a vector Wiener process) is a scalar function $V: D\to \mathbb{R}$ defined on a region $D\subset\mathbb{R}^n$ that is continuous,
positive-definite, twice continuously differentiable throughout $D,$ except possibly for the origin $x=0,$ and satisfies the condition $LV\leq 0$ (or $LV<0$), where $L$ is the differential generator associated with the stochastic system.

For the system (2.2) the operator $L$ have the form
\begin{equation}\label{eqn-2-4}
\begin{array}{ll}
\begin{aligned}
\begin{split}
 L=(ax+by)\frac{\partial}{\partial
x}+(cx+my)\frac{\partial}{\partial
y}+\frac{1}{2}\Big[(ex+fy)^2\frac{\partial^2}{\partial
x^2}+\\ +2(ex+my)(gx+hy)\frac{\partial^2}{\partial x \partial y}+(gx
+hy)^2 \frac{\partial^2}{\partial y^2}\Big].\\
\end{split}
\end{aligned}
  \end{array}
  \nonumber
\end{equation}

In the deterministic case ($e=f=g=h=0$), the operator $L$ coincides with the operator $\dot{V}$ of "the derivative of $V$ with respect to the system".

Now we turn to the formulations of the theorems on stochastic stability of systems (2.2) and (2.2), (2.3). Proofs  of  the proposed theorems are based on the construction of appropriate Lyapunov functions in quadratic forms and applying general stability theorems of Kushner  [16, p.~39] and Khasminskii [17, pp.~152, 185].

In the absence of "fluctuating" \, terms ($e=f=g=h=0$)
the stability conditions obtained  below turn into the necessary and sufficient Routh-Hurwitz conditions for the deterministic part of system (2.2).

\textbf{1.} One of the coefficients of the "random" \, part of the system (2.2) is nonzero and the others are zero.


\setcounter{theorem}{0}

\begin{theorem} \label{result1}
Let, in the system (2.2), $e\neq 0,$ $f=g=h=0.$ Then the zero solution ($x(t)\equiv 0,$ $y(t)\equiv 0$) of the system (2.2) is stable in probability if the following inequalities hold

 1) $a+m<0,$  $am-bc>0,$  $b\neq 0,$

 2) $e^2<[2(a+m)(bc-am)]/[m^2+(am-bc)].$

In addition, any solution ($x(t),  y(t)$) of system (2.2) possesses the property: $(x(t),  y(t))\to \{(x, y): x=0\}$ as $t\to\infty$ with probability 1.
\end{theorem}

\begin{theorem} \label{result2}
{\it Let $e\neq 0,$ $f=g=h=0.$ Then the zero solution ($x(t)\equiv 0,$ $y(t)\equiv 0$) of the system (2.2), (2.3) is stable in probability if at least one of the following conditions holds:}

 1) $m>0,$ $a+m<0,$ $am-bc>0,$ $b\neq 0,$
$\dfrac{e^2}{2}<min\{-(a+m),\lambda_0\},$ \textit{where $\lambda_0$  is a positive root of the quadratic equation $Q(\lambda)=0,$}
$$
Q(\lambda)=2m\lambda^2+(2m^2+3am-2bc)\lambda-(a+m)(bc-am);
$$
\textit{ or }

2) $m<0,$ $a+m<0,$ $am-bc>0,$ $b\neq 0,$
$\frac{e^2}{2}<min\{-(a+m), \dfrac{bc-am}{m}, \lambda_0\},$
\textit{where $\lambda_0$ is the smallest positive root of the
equation $Q(\lambda)=0$ in the case
$$
B\geq 2\sqrt{AC},
$$
where $$A=2m, \; B=2m^2+3am-2bc, \; C=(a+m)(am-bc);$$ if
$0<B<2\sqrt{AC},$ then $\lambda_0$ is absent from hypothesis 2),}

3) $m=0, a<0, bc<0, e^2<-a.$ 

\textit{ In addition, any solution ($x(t), y(t)$) of the system } (2.2), (2.3) \textit{ possesses the property: } $(x(t), y(t))\to \{(x,y):x=0\}$ \textit{ as } $t\to
+\infty$ \textit{ with the probability} 1.
\end{theorem}


\begin{remark}
The inequalities with respect to parameter $e$ in the conditions 1)--3) of Theorem 1 give  upper bounds of the white noise intensity, up to which the perturbed system (2.2) remains stable in probability. The right-hand sides of these inequalities determine the lower bound for the bifurcation value of the white noise intensity.
\end{remark}

The Theorems 1 and 2 are proved using the following Lyapunov function
$$
V_{e}(x,y)=(am-bc)x^{2}+(mx-by)^{2}.
$$

Mean square exponential stability conditions of  zero solution of the systems (2.2) and (2.2), (2.3) in the case of $e\neq 0,$ $f=g=h=0$ are given by the following theorems.

\begin{theorem} \label{result3}
\textit{Let $e\neq 0,$ $f=g=h=0.$ Then for the zero solution ($x(t)=0, y(t)=0$) of  the system (2.2) to be exponentially stable in mean square it is necessary and sufficient that the following conditions should satisfy:}

1) $a+m<0, \; am-bc>0$;

$2a$) $e^2<\min\left\{\dfrac{2(a+m)(bc-am)}{m^{2}+(am-bc)}, \;
\dfrac{p+\sqrt{p^{2}+c^{2}q}}{c^{2}}\right\}$ 
\\ \textit{if} $c\neq 0,$ \textit{where}
$$
p=(a+m)(c^{2}+m^{2}+am-bc)-2c(ac+bm),
$$
$$
q=4(am-bc)\left[(a+m)^2+(b-c)^2\right];
$$

$2b$) $e^2<-2a$ $(a<0),$ \; \textit{if} $c=0.$ \\ 
\textit{In this case any solution $(x(t), y(t))$ of the system (2.2) has the property: $x(t)\to 0,$ $y(t\to 0$ as $t\to 0$ with probability} 1.
\end{theorem}

\begin{remark}
The inequalities $2a$) and $2b$) of theorem 3 give the bifurcation values of white noise intensity $e^2.$
\end{remark}

\begin{theorem} \label{result4}
\textit{Let} $e\neq 0,$ $c=f=g=h=0.$ \textit{Then the zero solution $(x(t)\equiv 0, y(t)\equiv 0)$ of system (2.2), (2.3) is exponentially stable in mean square if and only if} $a<0,$ $m<0,$ $e^2<-2a.$
\\ \textit{In this case any solution ($x(t), y(t)$) of the system} (2.2), (2.3) \textit{has the property:} $x(t)\to 0,$ $y(t)\to 0,$ \textit{as} $t\to +\infty$ \textit{with probability} 1.
\end{theorem}

\begin{theorem} \label{result5}
\textit{Let} $c\neq 0,$ $m=0,$ $e\neq 0,$ $f=g=h=0.$ \textit{Then the zero solution $(x(t)\equiv 0,
y(t)\equiv 0)$ of system (2.2), (2.3) is exponentially stable in mean square if and only if at least one of the following two conditions is satisfied:}

1) $a<0,$ $bc<0,$ $(-b)/c\geq 1,$ $e^2<-a$ \textit{or}

2) $a<0,$ $bc<0,$ 0<$(-b)/c< 1,$
$\dfrac{e^2}{2}<\dfrac{(-a)(b+c)}{b+3c}.$ \\ \textit{In this case any solution} $(x(t), y(t))\to (0,0)$ \textit{as}  $t\to +\infty$ \textit{with probability} 1.
\end{theorem}

\begin{remark}
The comparison between the stability conditions given by Theorems 1 and 2, and  Theorems 3 and 4 (or 5, in the hypotheses of theorems 4 and 5) shows that,  in general, the stability conditions  for the system (2.2) in Ito's  form are
wider than  the stability conditions  for this system in Stratonovich's one.
\end{remark}

In the case $f\neq 0,$ $e=g=h=0$ the conditions of stability in probability and mean square exponential stability of zero solution of the (2.2) in the Ito form coincide with the corresponding stability conditions of the system (2.2), (2.3) in the
Stratonovich form. The following theorems are valid.

\begin{theorem} \label{result6}
\textit{Let} $f\neq 0,$ $e=g=h=0.$
\textit{Then a sufficient condition for stability in probability of zero solution ($x(t)=0, y(t)=0$) of the system (2.2) is that the following inequalities should be satisfied:}

1) $a+m<0,$ $am-bc>0,$ $c\neq 0;$

2) $f^{2}<[2(a+m)(bc-am)]/c^{2}.$ \\ \textit{In this case any
solution} $(x(t),y(t))\to \{(x,y):y=0\}$ \textit{as $t\to\infty$ with probability} 1.
\end{theorem}

The statement of theorem 6 is proved using the Lyapunov function
$$
V_{f}(x,y)=(cx-ay)x^{2}+(am-bc)y^{2}.
$$

\begin{theorem} \label{result7}
\textit{Let} $f\neq 0,$ $e=g=h=0.$
\textit{Then a necessary and sufficient condition for mean square exponential stability of zero solution ($x(t)=0, y(t)=0$) of the system (2.2) is that the following  inequalities should be fulfilled:}

1) $a+m<0, \; am-bc>0$;

2) $f^{2}<\min\left\{\dfrac{2(a+m)(bc-am)}{c^{2}}, \;
\dfrac{-p+\sqrt{p^{2}+c^{2}q}}{c^{2}}\right\}$ \\ \textit{ in the case } $c\neq 0,$
\textit{where}
$$
p=(a+m)(c^{2}+m^{2}+am-bc)-2c(ac+bm),
$$
$$
q=4(am-bc)\left[(a+m)^2+(b-c)^2\right];
$$
\textit{If $c=0,$ then the condition} 2) \textit{is lacking.}

\textit{In addition, any solution $(x(t),y(t))\to (0,0)$ as $t\to \infty$ with probability} 1.
\end{theorem}

The Theorem 7 is proved using the Lyapunov function
\begin{displaymath}
\begin{array}{ll}
V_{f}(x,y)=[c^{2}+m^{2}+(am-bc)]x^{2}+2(ac+bm-\\ -c{f^{2}}/{2})xy+[a^{2}+b^{2}+(am-bc)-(a+m){f^{2}}/{2}]y^{2}.
\end{array}
\end{displaymath}

The cases $g\neq 0;$ $e=f=h=0$ and $h\neq 0;$ $e=f=g=0$ are reduced to the cases discussed above, by replacing $x\to y,$ $y\to x,$ $a\to m,$ $b\to c,$ $c\to b,$ $m\to a$ and $g\to f$ in the first case and $h\to e$ in the second case, respectively.

\textbf{2.} Two coefficients of the "random" part of the system (2.2) are different from zero and the rest are zero. To avoid bulky expressions, we limit ourselves to the cases when these two non-zero coefficients from $e, f, g, h$ are  equal to each other.

The following statement takes place.

\begin{theorem} \label{result8}
\textit{Let $e=f\neq 0,$ $g=h=0.$ Then a sufficient condition for the zero solution ($x(t)=0, y(t)=0$) of the system (2.2) to be stable in probability is that the following inequalities should be satisfied:}

1) $a+m<0,$ $am-bc>0,$ $c\neq 0;$

2) $e^2< min \{[2(a+m)(bc-am)]/[(am-bc)+(c-m)^2];$

$(m-a-2c)+\sqrt{(a+m)^2+4c(a-m+c-b)}\}.$ \\ \textit{In this case any solution $(x(t),y(t))$ of system (2.2) $(x(t),y(t))\to \{(x,y):y=0\}$ as $t\to +\infty$ with probability} 1.
\end{theorem}

To prove Theorem 8 we use the Lyapunov function in the form
$$
V(x,y)=c^2x^2-2c(2a+c^2)xy+\Big[a^2+(am-bc)+\frac{e^2}{2}(a+m-2c)\Big]y^2.
$$

The following theorem gives a necessary and sufficient stability conditions in the mean square.

\begin{theorem} \label{result9}
{\it Let $e=f\neq 0,$ $g=h=0.$ Then the zero solution $(x(t)\equiv 0, y (t)\equiv 0)$ of the system (2.2) is exponentially stable in the mean square if and only if the the following inequalities are satisfied:

1) $a+m<0,$ $am-bc>0;$

2a) $e^2<\min \{[2(a+m)(bc-am)]/[(am-bc)+(c-m)^2],$ \\
 $(p+\sqrt{p^2+4m^2q})/(2m^2)\}$  if $m\neq 0,$  where
$$
p=(b-c)(b^2+c^2+am-bc)-2m(ac+bm),
$$
$$
q=(am-bc)[(a+m)^2+(b-c)^2];
$$

$2b$) $e^{2}<\min \{(2ab)/(c-b);$ $((-b)
[a^2+(b-c)^2])/((c-b)^2)\}$ \; if $m=0$ and $b<0,$ $c>0$ \\ (if $b>0,$ $c<0,$ then the condition $2b$) is missing). 

In addition, any solution $(x(t),y(t))$ of system (2.2)  has the property $x(t)\to 0,$ $y(t)\to 0$ as $t\to +\infty$ with probability} 1.
\end{theorem}

\begin{theorem} \label{result10}
\textit{Let} $e=f\neq 0,$ $g=h=0.$ \textit{Assume that $m=0,$ $a<0,$ $b<0,$ $c>0.$ Then for the  zero solution ($x(t)=0, y(t)=0$) of  the system (2.2), (2.3) to be exponentially stable in mean square it is necessary and sufficient that}
$$
\frac{e^2}{2}<min\Big\{\frac{ab}{c-2b},\frac{ab+(c-b)^2-(c-b)\sqrt{c^2+2b(a-c)}}{-b}\Big\}.
$$
\textit{In this case for any solution $(x(t),y(t))$ we have: $x(t)\to 0,$ $y(t)\to 0$ as $t\to +\infty$ with probability} 1.
\end{theorem}

\begin{theorem} \label{result11}
\textit{Let} $e=f\neq 0,$ $g=h=0.$ \textit{Suppose that $b=c=m,$ $a<b<0.$ Then in order that for the zero solution of the system} (2.2), (2.3) \textit{to be exponentially stable in the mean square, it is necessary and sufficient that the inequality should be satisfied}
$$
\frac{e^2}{2}<min\Big\{-\frac{a+b}{3},\frac{-3a-\sqrt{a^2+8b^2}}{4}\Big\}.
$$
\textit{In this case any solution $(x(t),y(t))\to (0,0)$ as $t\to +\infty$ with probability} 1.
\end{theorem}

The Theorems 9--11 are proved by using the Lyapunov function in the form 
$$
\begin{array}{ccc}
  V(x,y)=(c^{2}+m^{2}+am-bc)x^{2}+2(me^{2}+ac+bm)+ \\
 +\left(a^{2}+b^{2}+am-bc+(b-c)e^{2}\right)y^{2}.
\end{array}%
$$

\begin{remark}
The right-hand sides of the inequalities in conditions 2$a$) and 2$b$) of Theorem 9 and the inequalities in Theorems 10 and 11 give a bifurcation value of the white noise intensity.
\end{remark}

The following theorems provide sufficient conditions for stability in probability.

\begin{theorem} \label{result12}
{\it Let $e=g\neq 0,$ $f=h=0.$ Then sufficient conditions for stability in probability of the trivial solution $(x(t)\equiv 0, y(t)\equiv 0)$ of the system (2.2) are that the following inequalities should hold:

1) $a+m<0,$ $am-bc>0$ \; $b\neq 0;$

2) $e^2<[2(a+m)(bc-am)]/[(am-bc)+(b-m)^{2}].$ \\ In this case any solution $(x(t),y(t))$ of the system (2.2) possesses the property $(x(t),y(t))\to \{(x,y):y=0\}$ as $t\to +\infty$ with probability} 1.
\end{theorem}

\begin{theorem} \label{result13}
\textit{Let} $e=g\neq 0,$ $f=h=0.$ \textit{Then the zero solution of the system} (2.2), (2.3) \textit{is stable in probability if}
\begin{equation}
\begin{array}{ll}
  1) \; m<0, \; b\neq 0, \\
  2) \; e^2/2<min\{(am-bc)/(-m); \lambda_0\},
\end{array}%
 \nonumber
\end{equation}
\textit{where $\lambda_0$ is the positive root of the quadratic equation}
$$
2m\lambda^2+[2(am-bc)+m(a+m)+(b-m)^2]\lambda-(a+m)(am-bc)=0.
$$
\textit{In addition, any solution $(x(t),y(t))\to (0,0)$ as $t\to +\infty$ with probability} 1.
\end{theorem}

To prove the Theorems 12, 13 the Lyapunov function
$$
\begin{array}{ccc}
  V(x,y)=(am-bc)x^2+(mx-by)^2
\end{array}%
$$
is used.

The following theorem considers the case $e=h\neq 0,$ $f=g=0.$

\begin{theorem} \label{result14}
{\it Let $e=h\neq 0,$ $f=g=0.$ Then the trivial solution $(x(t)\equiv 0, y(t)\equiv 0)$ of system (2.2) is stable in probability if:

1) $a+m<0,$ \; $am-bc>0;$ $c\neq 0;$

2) $e^2<\min \{p_1, p_2\},$ \; \\ where $p_1$ and $p_2$ are the smallest positive roots of the equations:
$$
\lambda^2+2(a+m)\lambda+4(am-bc)=0
$$ under the condition $(a+m)^2\geq 4(am-bc)$
and
$$
\lambda^3+3(a+m)\lambda^2+2[(a+m)^2+2(am-bc)]\lambda+4(a+m)(am-bc)=0,
$$
respectively.

(If $(a+m)^2<4(am-bc),$ then the first equation has no real roots and therefore $p_1$ is absent.) 

In addition, any solution $(x(t),y(t))\to \{(x,y):y=0\}$ as $t\to +\infty$ with probability 1.}
\end{theorem}

The Theorem 14 is proved by using the following Lyapunov function
$$
V(x,y)=2c^2x^2-2c(2a+e^2)xy+[2(a^2+am-bc)+e^4+(3a+m)e^2]y^2.
$$

The next theorem considers the case $f=g\neq 0,$ $e=h=0.$

\begin{theorem} \label{result15}
{\it Let $f=g\neq 0,$ $e=h=0.$ Then the sufficient conditions for the trivial solution $(x(t)\equiv 0, y(t)\equiv 0)$ of the system (2.2)  to be stable are that the following inequalities should be satisfied:

1) $a+m<0,$ $am-bc>0;$

$2a$) $m<0, \; f^2<\min \{p_1, p_2, p_3 \},$ \; \\ where
$p_1=[m^2+(am-bc)]/(-m),$ $p_2$ and $p_3$ are the smallest positive roots of the equations, respectively:
$$
\lambda^3+(a+m)\lambda^2-2(b^2+c^2+2am)\lambda+4(a+m)(bc-am)=0,
\eqno(A)
$$
\begin{displaymath}
\begin{array}{ll}
2m\lambda^3+[c^2+2m^2+2(am-bc)+2m(a+m)]\lambda^2+\\
+2[(a+m)(m^2+am-bc)-4bm(b+c)]\lambda-4b^2(am-bc)=0.
\end{array}\eqno(B)
\end{displaymath}

(The equations ($A$) and ($B$) can not possibly have positive roots: in this case $p_2$ and/or $p_3$ are missing.)

$2b$) $m>0,$ $f^2<\min\{p_2, p_3\}.$

 In addition, any solution $(x(t), y(t))\to\{(x,y):x=0\}$ as $t\to+\infty$ with probability 1.}
\end{theorem}

The Theorem 15 is proved by using the Lyapunov function:
\begin{displaymath}
\begin{array}{rr}
V(x,y)=2[m^2+(am-bc)+mf^2]x^2-2(2bm-cf^2)xy+ \nonumber \\ +[2b^2-(a+m)f^2-f^4]y^2.
\end{array}
\end{displaymath}

\begin{remark}
The right-hand sides of inequalities in hypotheses 2 and 2$a$), 2$b$) in Theorems 14 and 15, respectively, give the lower bounds of the  bifurcation value of the white noise intensity $f^2.$
\end{remark}

The remaining two cases $f=h\neq 0,$ $e=g=0$ and $g=h\neq 0,$ $e=f=0$ are reduced to the considered cases  $e=g\neq 0,$ $f=h=0$ and $e=f\neq 0,$ $g=h=0$ by replacing $x\to y,$ $y\to x;$ $a\to m,$ $b\to c,$ $c\to b$, $m\to a$ and $h\to e,$ $f\to g,$ $g:=e$
(in the first case) and $h\to e,$ $g\to f,$ $f:=e$ (in the second case).

\section{Example: The damped harmonic oscillator} 

Consider a harmonic oscillator with eigenfrequency $\omega$ subject to the action of a damping force proportional to velocity with coefficient $k.$  In many problems it seems natural to assume that $k$ and $\omega$ are merely the mean value of the damping coefficient and eigenfrequency, while their true values are stochastic processes with small correlation interval. Such a system can be described by stochastic equation
$$
\ddot{x}(t)+(k+\sigma_2 \dot{\xi}(t))\dot{x}(t)+(\omega^2+\sigma_1
\dot{\xi}(t))x(t)=0, \eqno(4.1)
$$
where $x(t)$ is a displacement of a material point from the equilibrium position at the time $t$; $\dot{\xi}(t)$ is a white noise of unit intensity, $\sigma_1^2$ and $\sigma_2^2$ are the intensities of white noises acting on the frequency and damping coefficient, respectively.

We interpret the equation (4.1) in two senses:  as an Ito and a Stratonovich stochastic equations.

Denoting $\dot{x}(t)=y(t),$ rewrite the equation (4.1) as a system of stochastic equations
\begin{equation}
\left\{ \begin{array}{ll}
 dx(t)=y(t)dt,\\
  dy(t)=-[\omega^2x(t)+ky(t)]dt-[\sigma_1 x(t)+\sigma_2 y(t)]d \xi(t),\\
\end{array} \right.\nonumber
\eqno(4.2)
\end{equation}
where $d \xi(t)$ is the stochastic differential of the Wiener process $\xi(t)$ in the Ito sense.

The equation (4.1) interpreted in the Stratonovich form of stochastic equations is equivalent  to  the system of Ito's form
\begin{equation}
\left\{ \begin{array}{lr}
 d x(t)=y(t)dt,\\
 d y(t)=-\Big[\Big(\omega^2-\dfrac{\sigma_1\sigma_2}{2}\Big)x(t)+(k-\dfrac{\sigma_2^2}{2})y(t)\Big]dt-\\ \;\;\;\;\;\;\;\;\;\;\;\; -[\sigma_1 x(t)+\sigma_2 y(t)]d \xi(t).\\
 \end{array} \right.\nonumber
\eqno(4.3)
\end{equation}

We formulate the conditions of exponential stability in the mean square for the systems (4.2) and (4.3).

\textbf{Proposition 1.} {\it The trivial solution} ($x(t)\equiv 0, y(t)\equiv 0$) {\it of the system} (4.2) {\it is exponentially stable in the mean square if and only if the following conditions are satisfied:}
\begin{equation}
\begin{array}{lll}
   a) \; \sigma_1^2<2k\omega^2 \textit{  in the case  } \sigma_1\neq 0, \, \sigma_2=0,  \\
   b) \; \sigma_2^2<2k \textit{  in the case  } \sigma_1=0, \, \sigma_2\neq 0 \\
   c) \; \sigma^2<\dfrac{2k\omega^2}{\omega^2+1} \textit{  in the case  } \sigma_1=\sigma_2=\sigma. \\
\end{array}%
\nonumber
\end{equation}

The proposition 1 is proved by using the Lyapunov functions:
\begin{equation}
\begin{array}{lll}
   a) \; V(x,y)=\dfrac{k^2+\omega^2(\omega^2+1)+k
{\sigma_1^2/2}}{2k\omega^2-\sigma_1^2}x^2+
\dfrac{2k+\sigma_1^2}{2k\omega^2-\sigma_1^2}xy+\\ +\dfrac{\omega^2+1}{2k\omega^2-\sigma_1^2}y^2,  \\
   b) \; V(x,y)=\left( \dfrac{\omega^2+1}{2k-\sigma_2^2}+\dfrac{k}{2\omega^
2}\right)x^2+\dfrac{1}{\omega^2}xy+\dfrac{\omega^2+1}{\omega^2(2k-\sigma_2^2)}y^2, \\
   c) \; V(x,y)=\dfrac{1}{2[\sigma^2(\omega^2+1)-2k\omega^2]} \Big\{
2[(\omega^2+1)\sigma^2-\omega^4-\omega^2-\\ -k^2]x^2-4kxy-2(\omega^2+1)y^2\Big\}. \\
\end{array}%
\nonumber
\end{equation}

From Proposition 1 and  Kushner's Theorem [16, p. 39] we obtain the following statement.

\textbf{Corollary.} {\it For each number $p>0$ and for any solution $(x(t), y(t))$ with the initial condition $x(0)=x_0,$ $y(0)=y_0,$ where $(x_0, y_0)\in \{(x,y):V(x,y)<p\},$ the following relations hold: $x(t)\to 0,$  $y(t)\to 0$ as $t\to +\infty$ with probability at least $1-(V(x_0)/p),$  where $V$ is the Lyapunov function corresponding to the cases} $a$), $b$), $c$), respectively.

\begin{remark}
Since the number $p$ is arbitrary we have $x(t)\to 0,$ $y(t)\to 0$ with probability 1.
\end{remark}

Using the Lyapunov functions similar to those  $a$)--$c$) above, we get the following

\textbf{Proposition 2.} \textit{The trivial solution of the system} (4.3) \textit{is exponentially stable in the mean square  if and only if:}
\begin{equation}
\begin{array}{lll}
   a) \; \sigma_1^2<2k\omega^2 \textit{ in the case  } \sigma_1\neq 0, \sigma_2=0,  \\
   b) \; \sigma_2^2<k \textit{ in the case  } \sigma_1=0, \sigma_2\neq 0, \\
   c) \; \sigma^2<\omega^2+\frac{k+1}{2}-\sqrt{\omega^4+(1-k)\omega^2+\frac{(k+1)^2}{4}} \\ \textit{ in the case  } \sigma_1=\sigma_2=\sigma\neq 0. \\
\end{array}%
\nonumber
\end{equation}
\textit{In addition for any solution $(x(t),y(t))$ of the system} (4.3) $x(t)\to 0,$ $y(t)\to 0$ \textit{as $t\to +\infty$ with probability} 1.

\begin{remark}
The right-hand sides of the inequalities $a$)--$c$) of the propositions 1 and 2 determine the exact values of the bifurcation values of the intensities $\sigma_1^2,$ $\sigma_2^2$ and $\sigma^2$ of the white noise.
\end{remark}

\begin{remark}
Comparing the inequalities $a$)--$c$) from propositions 1 and 2 we see that in general the stability conditions for equation (4.2) in the Ito form are wider than the corresponding stability conditions for the equation (4.3) in the Stratonovich form.
\end{remark}

\section{Conclusion} 

In the paper sufficient conditions for stability in probability,  necessary and sufficient conditions for exponential stability in the mean square of the zero solution of two-dimensional linear autonomous stochastic systems in the Ito and Stratonovich forms are obtained. The boundaries of the bifurcation values of the white noise intensity, which affects the parameters of the stochastic system  are established. As an example a damped  harmonic linear oscillator equation is considered in the case when both the natural frequency and friction coefficient are affected by white noise. Necessary and sufficient conditions for exponential stability of trivial solution in the mean square of the perturbed harmonic oscillator are given.

\end{document}